\documentclass[12pt]{article}
\usepackage{amsmath,amsfonts,amssymb,amsthm,epsfig}
\newcommand{\bb}{\mathbb}
\newcommand{\R}{\bb  R}
\newcommand{\C}{\bb C}

\newcommand{\M}{\mathcal M}
\newcommand{\A}{\mathcal A}
\newtheorem{theorem}{Theorem}

\newtheorem{proposition}{Proposition}

\newtheorem{example}{Example}
\begin{document}

\title{A divergent Teichm\"uller geodesic with 
       uniquely ergodic vertical foliation}  
\author{
   Yitwah Cheung\\
   Northwestern University\\ Evanston, Illinois\\
   email: \emph{yitwah@math.northwestern.edu}
\bigskip\\ and \bigskip\\
   Howard Masur\thanks{This research is partially 
   supported by NSF grant DMS0244472.}\\
   University of Illinois at Chicago\\ Chicago, Illinois\\
   email: \emph{masur@math.uic.edu}
}

\date{\today}

\maketitle
\begin{abstract}
\noindent  We construct an example of a quadratic differential whose vertical foliation is uniquely ergodic and such that the Teichm\"uller geodesic determined by the quadratic differential diverges in the moduli space of Riemann surfaces. 
\end{abstract}

\section{Introduction and Statement of theorem}
We let $S$ be a closed surface of genus $g\geq 2$ and $T_g$ the \emph{Teichm\"uller space} of $S$ equipped with the Teichm\"uller metric.  Let $\text{Map}=\text{Diff}^+(S)/\text{Diff}_0(S)$  denote the \emph{mapping class group} of $S$.  It acts 
on $T_g$ by isometries with quotient space $\M_g$,
 the \emph{moduli space} of surfaces of genus $g$.  It is well-known that $\M_g$ is not compact; one  leaves compact sets of $\M_g$ by finding curves whose lengths in the hyperbolic metric on the surfaces  
 goes to zero.
For any Riemann surface $X\in T_g$ a \emph{holomorphic quadratic differential}  $q=q(z)dz^2$ assigns to each  uniformizing parameter $z$  a holomorphic function $q(z)$ such that $q(z)dz^2$ is invariant under holomorphic change of coordinates. Away from the zeroes of $q$ there are natural coordinates so that $q(z)\equiv1$  so that  $q$ defines a metric $|dz^2|$ which is locally Euclidean away from the zeroes of $q$.  A \emph{saddle connection} $\gamma$ is a geodesic in this metric joining a pair of (not necessarily distinct) zeroes that has no zeroes in its interior. It is represented by a straight line in the natural Euclidean metric and determines a \emph{holonomy vector} with horizontal and vertical components which we denote by $\lambda(\gamma)$ and $v(\gamma)$ respectively.

Associated to $q$ are the \emph{horizontal} and \emph{vertical} trajectories.  These are the arcs along which $q(z)dz^2>0$ and $q(z)dz^2<0$. 
For each $t\in\R$ one can define a new quadratic differential $g_t(q)$ on a new Riemann surface $X_t$
by expanding along the \emph{horizontal} trajectories  of 
$q$ by a factor of $e^t$ and contracting along the vertical trajectories  by a factor of $e^{-t}$.  The family of $X_t$ form a \emph{Teichm\"uller geodesic}  through $X$ whose  projection to $\M_g$ defines a geodesic in $M_g$.  The underlying map $X\to X_t$ of Riemann surfaces is called a \emph{Teichm\"uller map}.
In the appropriate natural coordinates $z=x+\sqrt{-1}y$  on $X$ away from the zeroes of $q$, and natural coordinates $z_t=x_t+\sqrt{-1}y_t$ on $X_t$ away from the  zeroes of $g_t(q)$, the map is given by 
$$x_t=e^tx,\\ y_t=e^{-t}y.$$

One may study the dynamics of the family of vertical trajectories of $q$ which defines  the  \emph{vertical foliation} $F_q$. 
An important notion in topological dynamics is that of minimality. The foliation $F_q$ is \emph{minimal} if the full orbit of every leaf is dense.  It is a standard fact
\cite{St} that if there are no  saddle connections with zero horizontal holonomy,   then  the vertical foliation $F_q$ is minimal. 
  
For minimal foliations it is natural to study their ergodic behavior. A foliation is \emph{uniquely ergodic} if every leaf is uniformly distributed on the surface. Equivalently, the foliation is uniquely ergodic if there is a unique, up to scalar multiplication,   measure transverse to $F_q$  invariant under holonomy along leaves of $F_q$.  (To avoid trivialities, we assume these measures  are  supported on the complement of the singularities)   
An interesting phenomenon  found by \cite{S}, \cite{V1},
\cite{Ke},\cite{KN}   is the existence of minimal foliations that are not uniquely ergodic.  The last two  examples were in the context of interval exchange transformations, but a suspension of  an interval exchange transformation yields a 
(oriented) measured foliation.   

In  \cite{M} a connection was made between the dynamics of the foliation $F_q$ and the dynamics of the corresponding geodesic $X_t$ determined by $q$ in $\M_g$.  It was shown that if the foliation is minimal, but not uniquely ergodic, then the  geodesic \emph{diverges} in $\M_g$.   This means that it eventually leaves every compact set.  The question  arises if this condition is necessary; does  $X_t$ divergent imply $F_q$ is  nonuniquely ergodic?
We show the answer is negative.
\begin{theorem}
\label{th:main}
There exist quadratic differentials $q$ with uniquely ergodic $F_q$ such that the Teichmuller geodesic $X_t$ diverges in $\M_g$.
\end{theorem}

We will construct an example in genus two, although our  method will produce  
examples in any genus.   That such examples should exist was already suggested by Kerckhoff to the second author (oral communication).
We consider the  ``zippered rectangle'' construction of Veech \cite{V2}. 
 This  is a way of defining an Abelian differential on a Riemann surface 
such that the first return map of the vertical trajectories of the Abelian differential to a horizontal  interval is a given   interval exchange
transformation.   

   We then recall the notion of Rauzy induction  (\cite{R}.) This is a procedure which by taking  the first return map on a subinterval of the horizontal 
interval  gives a new interval exchange and also gives a new zippered rectangle. 
Associated to Rauzy induction is a graph of permutations of  interval 
exchanges.  The point of this construction in our context is that  
  successive applications of Rauzy induction gives a 
discrete set of points along the  Teichm\"uller geodesic defined by the 
zippered rectangles.

We then specialize to interval exchanges on four intervals and explicitly find the graph of permutations.   We produce an infinite  path in this graph that 
allows us to explicitly compute the lengths and heights of the zippered rectangles given by  Rauzy induction. 
The  interval exchange  
whose sequence   of Rauzy inductions give this path has the desired properties. Namely, it 
 is uniquely ergodic, and any Abelian differential formed by the zippered rectangle construction yields a divergent geodesic.


\section{Interval exchanges and zippered rectangles}
Let $\lambda\in\R_+^m$ be a vector of positive lengths and $\pi$ 
an irreducible permutation on $\{1,\ldots,m\}$ for some $m\geq2$.  
Recall irreducibility means $\pi\{1,\ldots,k\}\neq\{1,\ldots,k\}$ 
for $1\leq k<m$.  Associated to the pair $(\lambda,\pi)$ is an 
interval exchange map $T$ defined as follows.  
Let $I_j=[\beta_{j-1},\beta_j)$ where $\beta_0=0$ and 
$\beta_j=\lambda_1+\cdots+\lambda_j$ for $j=1,\ldots,m$.  
Then $T$ is the map defined on $I=\cup I_j$ by the formula: 
for $x\in I_j$
$$T(x)=x-\sum_{i<j}\lambda_i+\sum_{\pi i<\pi j}\lambda_{\pi j}.$$  

We recall a construction in \cite{V2} of  ``zippered rectangles'' to 
suspend $T$.  Let $h,a\in\R^m$ be vectors satisfying the inequalities 
\begin{align}
         h_j &> 0 &(1\leq j \leq m)\\
     0 < a_j &\leq \min(h_j,h_{j+1}) &(1\leq j<m, j\neq\pi^{-1}m)\\
   0 < a_{\pi^{-1}m} &\leq h_{\pi^{-1}m+1} &\\
   -h_{\pi^{-1}m} \leq a_m &\leq h_m &
\intertext{and the system of linear equations}
   h_j - a_j &= h_{\sigma j+1} - a_{\sigma j} &(0\leq j\leq m)
\end{align}
where $a_0=h_0=h_{m+1}=0$ by convention and $\sigma$ is the permutation 
on $\{0,1,\ldots,m\}$ defined by 
\begin{align}
   \sigma j = \left\{\begin{tabular}{ll}
              $\pi^{-1}1-1$ & $j=0$\\
              $m$           & $j=\pi^{-1}m$\\
              $\pi^{-1}(\pi j+1)-1$ & otherwise
              \end{tabular}\right..
\end{align}

The surface $M(\pi,\lambda,h,a)$ is obtained by glueing together 
$m$ rectangles along their boundaries in the following manner.  
Let $R_j$ be the rectangle in $\C$ with base $I_j$ on the real 
axis and height $h_j$.  There are three cases depending on the 
sign of $a_m$ and in each case there will be three sets of 
identifications, all of which are the form $z\to z+c$.  

First, consider the case $a_m=0$.  The first set of identifications is 
\begin{enumerate}
  \item[(1)] the top of $R_j$ is glued to the interval $T I_j$ 
             at the base for $j=1,\ldots,m$.  
\end{enumerate}
To describe the remaining identifications we shall use the notation 
$$R^+_j[a,b]\sim R^-_k[c,d]$$ to mean the right side of $R_j$ between 
heights $a$ and $b$ is glued to the left side of $R_k$ between heights 
$c$ and $d$.  (For this to be well-defined, we must have $b-a=d-c$.)  
The remaining identifications in the case $a_m=0$ are 
\begin{enumerate}
  \item[(2)] $R^+_j[0,a_j]\sim R^-_{j+1}[0,a_j]$ for $j=1,\ldots,m-1$, and 
  \item[(3)] $R^+_j[a_j,h_j]\sim R^-_{j+1}[a_{\sigma j},h_{\sigma j+1}]$ 
             for $j=1,\ldots,m$.  
\end{enumerate}
Now consider the case $a_m>0$.  Let $j=\pi^{-1}m$.  All identifications 
remain the same except the $j$th in (2), which is replaced with 
\begin{align*}
R^+_{\pi^{-1}m}[0,h_j]&\sim R^-_{\pi^{-1}m+1}[0,h_j] \quad\text{and}\\
R^+_m[0,a_m]&\sim R^-_{\pi^{-1}m+1}[h_j,a_j].  
\end{align*}
Similarly, in the case $a_m<0$ all identifications remain the same as in 
the first case except the $m$th in (3), which is replaced with 
\begin{align*}
R^+_{\pi^{-1}m}[a_{\pi^{-1}m},h_{\pi^{-1}m}]&\sim 
R^-_{\sigma m+1}[a_{\sigma m},h_{\sigma m+1}-h_m]\quad\text{and}\\ 
R^+_m[0,h_m]&\sim R^-_{\sigma m+1}[h_{\sigma m+1}-h_m,h_{\sigma m+1}].  
\end{align*}

The collection of glued rectangles $M=M(\pi,\lambda,h,a)$ is called 
the \emph{zippered rectangle} associated to $(\pi,\lambda,h,a)$.  
Since the glueing maps are of the form $z\to z+c$, $M$ carries the 
structure of a Riemann surface and the $1$-form $dz$ induces an Abelian 
differential $\omega=\omega(\pi,\lambda,h,a)$ on $M$.  The interval 
exchange $T$ is the first return map to $I$ under the flow in the 
vertical direction generated by the vector field $\partial/\partial y$.  

Note if each cycle of $\sigma$ contains at least two elements in 
$\{1,...,m-1\}$, then $\omega$ has a zero at $(0,0)$ and 
at $\lambda_1+\ldots+\lambda_j+ia_j$ for $j=1,\ldots,m$.  
In particular, each $R_j$ has at most one zero on each of its vertical 
sides.  More specifically, the left side always contains a zero while 
the right side contains a zero except when $j=m$, $a_m<0$ or when 
$j=\pi^{-1}m$, $a_m>0$.  

\begin{example}
Let $m=4$, $\pi j=5-j$ and $\lambda\in\R_+^4$.  Then $(\pi,\lambda,h,a)$ 
is a zippered rectangle with $h=(1,3,3,2)$, $a=(2,2,2,1)$ or with 
$h=(2,3,3,1)$ and $a=(1,1,1,-1)$, where $\sigma j=j-2 \pmod 5$.  
See Figure 1.  
\end{example}

\begin{figure}
\begin{center}
\includegraphics{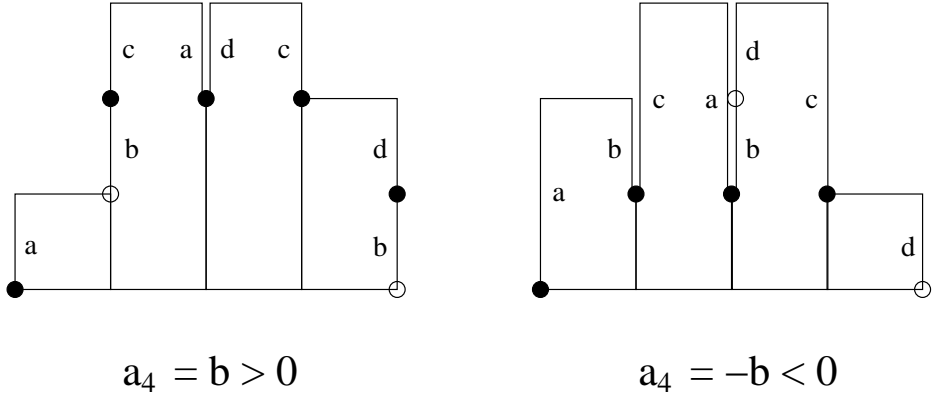}\label{F:Zipper}
\caption{Zippered rectangles in Example~1.}
\end{center}
\end{figure}

\section{Rauzy induction}
Let $\A$ denote an alphabet with $m\geq 2$ elements.  
A \emph{marked interval exchange} is a pair $(T,\nu_0)$ 
where $T$ is an interval exchange on $m$ intervals and 
$\nu_0$ a bijection from $\{1,\dots,m\}$ to $\A$.  
We think of the intervals of $T$ as being marked from 
left to right with names $\nu_0(1),\ldots,\nu_0(m)$.  
If $T$ is a $(\lambda,\pi)$ interval exchange, then the 
names of the image intervals from left to right are given by 
$\nu_1(1),\ldots,\nu_1(m)$ where $$\nu_1=\nu_0\circ\pi^{-1}.$$  

An alternative way to represent a marked interval exchange is 
via a triple $(\lambda,\nu_0,\nu_1)$ where $\lambda:\A\to\R_+$ 
is the function given by 
$$\lambda(\nu_0(j))=\lambda_j\qquad\text{for}\;j=1,\ldots,m.$$  
The function $\lambda$ and the vector $\lambda\in\R_+^m$ are 
represented by the same symbol so that $\lambda(\alpha)$ is the 
length of the interval marked $\alpha$, while $\lambda_j$ is the 
length of the $j$th interval from left to right.  We shall also 
refer to the pair $(\nu_0,\nu_1)$ as a \emph{marked permutation}.  
In what follows, we fix $\A=\{1,\ldots,m\}$.  

We define
\emph{Rauzy induction} on the space of marked interval exchange transformations. First for $\lambda$ a vector let $|\lambda|$ denote the sum of the lengths of the components.
Now given a  marked interval exchange $(\lambda,\nu_0,\nu_1)$, 
consider the  first return map 
found by inducing on the longer of the two subintervals 
$[0,|\lambda|-\lambda(\nu_i(m)))$ where $i\in{0,1}$.   
If $$\text{Case (a)}\qquad\lambda(\nu_1(m))>\lambda(\nu_0(m))$$ 
then the interval marked $\nu_1(m)$ is shortened.   The new marked permutation is $(\nu_0',\nu_1')$  where  $\nu'_1=\nu_1$ and 
$\nu'_0$ is the ordering obtained from $\nu_0$ by inserting $\nu_0(m)$ 
is the position immediately after $\nu_1(m)$ and moving forward one 
place the names of the intervals appearing after $\nu_1(m)$.   

Likewise, if $$\text{Case (b)}\qquad\lambda(\nu_0(m))>\lambda(\nu_1(m))$$ 
then the interval marked $\nu_0(m)$ is shortened, $\nu'_0=\nu_0$ 
and $\nu'_1$ is the ordering obtained from $\nu_1$ by inserting 
$\nu_1(m)$ in the position immediately after $\nu_0(m)$ and moving 
forward one place the names of the intervals appearing after $\nu_0(m)$.    
We ignore the case $\lambda(\nu_0(m))=\lambda(\nu_1(m))$. 

The result of Rauzy induction is a new marked interval exchange 
$(\lambda',\nu'_0,\nu'_1)$.  In the first case the lengths are 
related by 
\begin{equation}\label{E:Lengths}
\lambda=A\lambda',
\end{equation}
where $A$ is the elementary matrix whose only nonzero off-diagonal 
entry is a $1$ in the $(\nu_1(m),\nu_0(m))$ position.  In the second 
case the matrix is the transpose of the above matrix so that it has 
a $1$ in the $(\nu_0(m),\nu_1(m))$ position.  We can rephrase  the 
statement about the new marked permutation as follows. 
The new marked permutation 
is $a(\nu_0,\nu_1)$ in the first case and $b(\nu_0,\nu_1)$ in the second 
where $a,b$ are the bijections on the set of marked permutations defined by 
\begin{align}
  a(\nu_0,\nu_1) &= (\nu_0\circ c_k,\nu_1)
                    \qquad k=\nu_1^{-1}(m)\\
  b(\nu_0,\nu_1) &= (\nu_0,\nu_1\circ c_l)
                    \qquad l=\nu_0^{-1}(m)
\intertext{and}
  c_k(i) &= \left\{\begin{tabular}{ll}
            $i$ & $i\leq k$\\
            $m$ & $i=k+1$\\
            $i-1$ & $k+1<i\leq m$
            \end{tabular}\right..
\end{align}
(The context will make it clear whether $a$ is an operation on 
the set of marked permutations or a vector of lengths in $\R^m$.)  

By repeating this procedure, we obtain a sequence of elementary 
matrices $(A_k)_{k\geq1}$ and length vectors $(\lambda_k)_{k\geq1}$ 
satisfying $$\lambda=A_1\cdots A_k\lambda_k.$$  
The sequence $A_1\cdots A_k\cdots$ is called the \emph{expansion} 
of $\lambda$.  The corresponding sequence of marked permutations 
is a path in the extended Rauzy diagram.  This  is a directed graph 
whose vertex set is the extended Rauzy class $R(\nu_0,\nu_1)$ 
consisting of all marked permutations obtainable by applying the 
operations $a$ and $b$ to $(\nu_0,\nu_1)$.  There is a directed 
edge from $x$ to $y$ if and only if either $y=ax$ or $y=bx$.  The  
edges in the Rauzy diagram are labeled by $a$ or $b$ as the case may be.  
A directed path starting at $(\nu_0,\nu_1)$ is uniquely represented 
by a word in the alphabet $\{a,b\}$.  If $w$ is a word in $\{a,b\}$, 
let $\left<x;w\right>$ denote the path starting at $x$ obtained by 
following the letters of $w$ from left to right.  Let $[x;w]$ denote 
the corresponding product of elementary matrices.  If the path ends 
at $y$ and $w'$ is another word, then $[x;ww']=[x;w][y;w']$ where 
$ww'$ denotes the concatenation of the words $w$ and $w'$.  

\subsection{Rauzy induction with heights}
Rauzy induction also gives a transformation on the space of zippered 
rectangles.    Let 
$M(\pi,\lambda,h,a)$ be a zippered rectangle.
Let 
\begin{equation}\label{E:Heights}
h'=A^th
\end{equation}
and define, if $\lambda_m<\lambda_{\pi^{-1}m}$ (edge $a$)
\begin{equation}
   a'_j = \left\{\begin{tabular}{lr}
             $a_j$ & $j<\pi^{-1}m$\\
             $h_{\pi^{-1}m} + a_{m-1}$ & $j=\pi^{-1}m$\\
             $a_{j-1}$ & $\pi^{-1}m<j\leq m$
          \end{tabular}\right.
\end{equation}
and if $\lambda_m>\lambda_{\pi^{-1}m}$ (edge $b$)
\begin{equation}
   a'_j = \left\{\begin{tabular}{lr}
             $a_j$ & $j<m$\\
             $-(h_{\pi^{-1}m}-a_{\pi^{-1}m-1})$ & $j=m$
          \end{tabular}\right..
\end{equation}
Then the Riemann surface $M(\pi,\lambda,h,a)$ is biholomorphically 
equivalent to $M(\pi',\lambda',h',a')$ and the Abelian differential  
$\omega(\pi,\lambda,h,a)$ is equivalent to the Abelian differential 
$\omega(\pi',\lambda',h',a')$.  
In the first case, the last rectangle is stacked on top (as far 
right as possible) of a piece of the one that goes to the last.  
In the second case a piece of the last rectangle is stacked on top 
of the one that goes to the last.  See Figure 2.  
In particular, the last rectangle of the new set of zippered 
rectangles has a zero on its right side if the corresponding 
edge is labeled with an 'a' (the first case).  

\begin{figure}
\begin{center}
\includegraphics{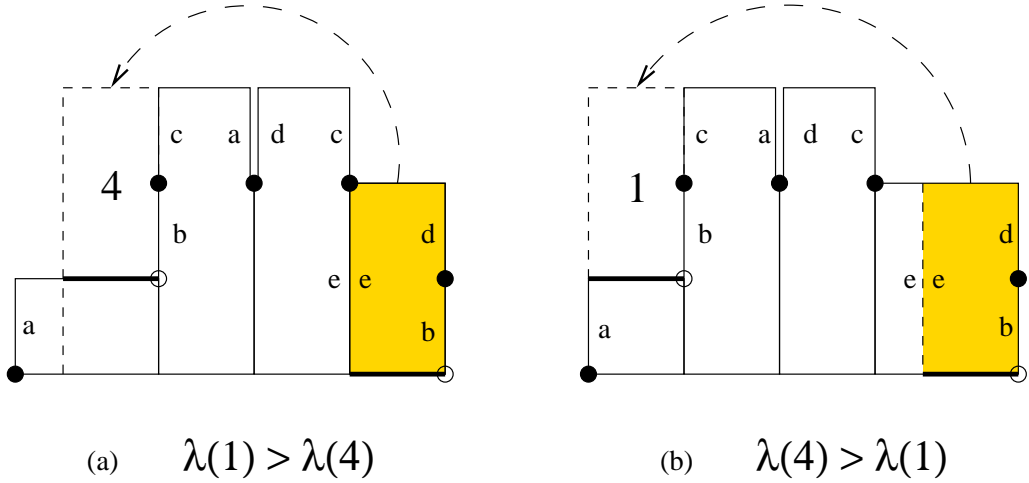}\label{F:Stacking}
\caption{Rauzy induction with $(\pi,\lambda,h,a)$ as in Example~1.  
The solid line indicates where a piece of the last rectangle is 
stacked on top of the first.  Note that the name of the rectangle 
with a new height is $4$ in the first case and $1$ in the second.}
\end{center}
\end{figure}

\section{Interval exchanges with four intervals}
For the rest of the paper we will consider interval exchanges on four 
intervals.  Let $\pi_0=(\nu_0,\nu_1)$ where $\nu_0(j)=j$ and $\nu_1(j)=5-j$.  
The extended Rauzy diagram $R(\pi_0)$ is shown in Figure 3.  

\begin{figure}
\begin{center}
\includegraphics{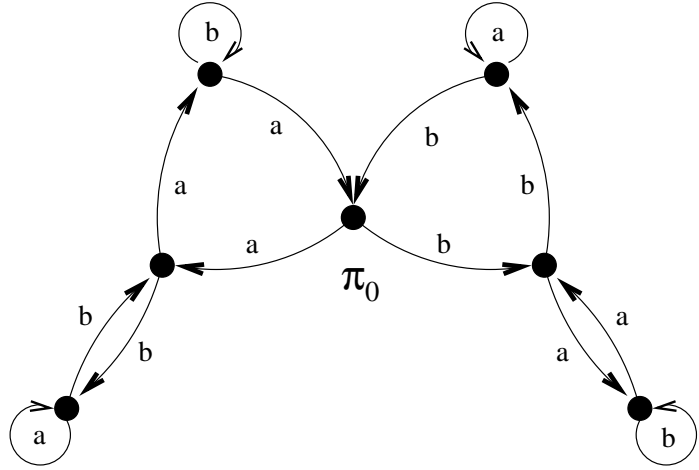}\label{F:Rauzy}
\caption{Rauzy diagram for the marked permutation $\pi_0$.}
\end{center}
\end{figure}

Now for each $n$ we form a pair of paths $$\left<\pi_0;bab^nab^2\right>$$ 
and $$\left<\pi_0,aba^nba^2\right>$$ that begin and end at $\pi_0$.  
We form the corresponding product of elementary matrices 
$A_n=[\pi_0;bab^nab^2]$ and  $B_n=[\pi_0;aba^nba^2]$ and set 
$C_n=A_nB_n$.  One computes $C_n$ as follows. At each stage on the path, the marked permutation $(\nu_0',\nu_1')$ is determined from the previous marked permutation $(\nu_0,\nu_1)$ by (8) and (9).  Each marked permutation determines an elementary  matrix as in (\ref{E:Lengths}) and the discussion that follows. The matrix $C_n$ is the product of these elementary matrices. It is then easily verified that the entries of $C_n$ are polynomials in $n$ 
whose leading terms are  
$$\begin{bmatrix}
  1 & 1 & 1 & 1\\
  1 & 2n & 2n & 1\\
  n &n^2&n^2& 0\\
  2 & n & n & 2\\
  \end{bmatrix}.$$

For each $n$, set  $$H_n=\prod_{j=1}^n C_j$$
For any 
$\lambda\in H_n\R_+^4$, let $\lambda_n$ be such that $$\lambda=H_n\lambda_n.$$  Then successive applications of Rauzy induction on the interval exchange $(\lambda,\pi_0)$ produces the
interval exchange $(\lambda_n,\pi_0)$.   
Now let $\Sigma$ be the standard simplex in $\R^4$ and regard each 
$C_j$ as a projective linear transformation from $\Sigma$ to itself. 
\begin{proposition}
\label{prop:path:uniquely}
The sets $\Sigma_n=H_n\Sigma$ are a decreasing sequence with infinite 
intersection a single point denoted by $\lambda_0$. The interval exchange $(\lambda_0,\pi_0)$ is uniquely ergodic.  
\end{proposition}
 In order to prove the proposition we use the Hilbert metric on $\Sigma$. We will prove a general proposition stronger than
what is needed. 
For any $m\times m$ matrix $L$ with all positive entries $L_{ij}$ set
 $$\delta(L)=\min_{i,j,r,s}
\frac{L_{ir}L_{js}}{L_{is}L_{jr}}.$$
\begin{proposition}
\label{prop:Hilbert}
Let $\Sigma$ be the open standard simplex in $\R^m$. Let $A_j:\Sigma\to\Sigma$ be a sequence of projective linear transformation defined by positive matrices with 
$$\prod_{j=1}^\infty(1-\delta(A_j))=0.$$  Then for any sequence $U_j$
of matrices with nonnegative entries,
$$\cap_{n=1}^\infty  A_1U_1A_2U_2\ldots A_nU_n\Sigma$$
is a single point. 
\end{proposition}
\begin{proof}
For $x,y\in\Sigma$  define
$$\Gamma(x,y)=\inf_{x_i,y_j\neq 0} \frac{x_jy_i}{x_iy_j}$$
and the Hilbert  metric
$$d(x,y)= -\log \Gamma(x,y).$$

Now for any matrix $U$ with
nonnegative entries, 
$$d(Ux,Uy)\leq d(x,y)$$ so it is enough to show that 
 $$d(Lx,Ly)\leq (1-\delta(L))d(x,y).$$
Since the metric is a path metric, it is enough to show this inequality for  $\epsilon:=d(x,y)$ arbitrarily small. By Lemma 15.2 of \cite{F}  $$\Gamma(Lx,Ly)\geq
\frac{\delta+\Gamma(x,y)}{1+\delta\Gamma(x,y)}$$
Writing $\Gamma(x,y)=e^{-\epsilon}=1-\epsilon+o(\epsilon)$,  we see that the right hand side of the above expression is of the form
$$\frac{\delta+1-\epsilon+o(\epsilon)}{1+\delta-\delta\epsilon+o(\epsilon)}.$$
and so  
$$\Gamma(Lx,Ly)\geq1-\frac{1-\delta}{1+\delta}\,\epsilon + o(\epsilon)$$
which implies that 
$$d(Lx,Ly)\leq 
\frac{(1-\delta)}{(1+\delta)}d(x,y)+o(d(x,y))$$ which gives the desired estimate.
\end{proof}

We now turn to the proof of Proposition~\ref{prop:path:uniquely}. 
From the form of the leading terms of the matrix $C_n$ we compute that the entries of the matrix $C_nC_{n+1}$ are all positive polynomials in $n$ whose leading terms are 
$$\begin{bmatrix}
  n & n^2 & n^2 & 4\\
  2n^2 & 2n^3 & 2n^3 & 2n\\
  n^3 &n^4&n^4& n^2\\
  n^2 & n^3 & n^3 & n\\
  \end{bmatrix}$$

We see immediately that there is an integer $p=p(i,j)\in [-2,2]$   such that the entries of  $C_nC_{n+1}$ satisfy 
$$1/4<n^p\frac{a_{i,r}}{a_{j,r}}<4$$ for any $r$, 
so that $$\delta(C_nC_{n+1})=1/16+O(1/n).$$ 
Proposition~\ref{prop:Hilbert} then implies that  
$$\text{diam} (H_n(\Sigma))\to 0$$ as $n\to\infty$.  Now it is a standard fact \cite{Ke} that if $\lambda_0=\cap H_n \Sigma$
 is minimal, it is uniquely ergodic.

We now show that $(\lambda_0,\pi_0)$  is minimal. 
 Fix a unit area  Abelian differential $\omega_0=\omega_0(\pi_0,\lambda_0,h_0,a_0)$.  
(This means that  the first return map of the vertical leaves to a horizontal 
segment defines the interval exchange $(\lambda_0,\pi_0)$ and $\omega_0$  defines a zippered rectangle with
height vector  $h_0=(h_0(1),\ldots,h_0(4))$.)

If $(\lambda_0,\pi_0)$ is not minimal then $\omega_0$ has a vertical saddle connection $\gamma$.  Choose $n$ large enough so that $\gamma$ does not intersect $\lambda_n$ in its interior. This means that there is some rectangle in the "zippered rectangles" determined by $\lambda_n$ that has $\gamma$ as part of its vertical side.  In that case there would be a pair of zeros on a side which contradicts the zippered rectangle construction.  

\section{Proof of Theorem~\ref{th:main}}
We show that  $\omega_0=\omega_0(\pi_0,\lambda_0,h_0,a_0)$ 
determines a divergent Teichm\"uller geodesic $X_t$. 
To do that we need to show that for all sufficiently large times $t$ along the geodesic $X_t$   there is a simple closed curve $\gamma_t$ which is short in the flat metric defined by  $g_t(\omega_0)$.  

The word $bab^nab^2aba^nba^2$ that is used to form the
matrix $C_n=A_nB_n$ ends with  $a$.  Recall that this means that 
the zippered rectangle for each interval exchange    $(\lambda_n,\pi_0)$
 has a zero on the right side of the fourth rectangle. 
Thus there is a saddle connection $\gamma_n$ such that the  horizontal component $\lambda(\gamma_n)$ of its holonomy is equal to $\lambda_n(4)$, and such that the vertical component $v(\gamma_n)$ of its holonomy is at most $h_n(4)$.  We will show that this  saddle connection 
 becomes short in some  time interval and that these time intervals cover 
$(T,\infty)$ for large enough $T$.  Since $\omega_0$ has only one singularity, 
this saddle connection will always be a closed curve.

We first give a bound for $\lambda_n(4)h_n(4)$.
By (\ref{E:Lengths}) the length vectors satisfy   
 the equation $$\lambda_n= C_{n+1}C_{n+2}\lambda_{n+2}$$
and by (\ref{E:Heights}) the height vector satisfies 
\begin{equation}
\label{eq:heights}h_n=(C_{n-1}C_n)^t h_{n-2}.
\end{equation}
The form of the leading terms of these matrices insures that for $i\neq 3$, 
\begin{equation}
\label{eq:three:wins}
\frac{\lambda_n(i)}{\lambda_n(3)}=O(1/n)
\end{equation}
 and 
\begin{equation}
\label{eq:three:doesnot lose}
\frac{h_n(i)}{h_n(3)}=O(1).  
\end{equation}

Since $$\sum_{i=1}^4 \lambda_n(i)h_n(i)=1,$$
the above estimates show that 
\begin{equation}
\label{eq:three:dominates}
\lambda_n(3)h_n(3)=1+O(1/n).
\end{equation}
Moreover the form of the leading terms of the matrices also gives 
$$\frac{h_n(4)}{h_n(3)}=O(1/n^2)$$ 
which together with (\ref{eq:three:wins}) and (\ref{eq:three:dominates}) gives 
\begin{equation}
\label{eq:saddle}
\lambda_n(4)h_n(4)=O(1/n^3).
\end{equation}

Now define  $s_n$ by 
\begin{equation}
\label{eq:s}
e^{-s_n}v(\gamma_n)=\frac{1}{\log n}
\end{equation}
 and $t_n$ by 
\begin{equation}
\label{eq:t}
e^{t_n}\lambda(\gamma_n)=\frac{1}{\log n}.
\end{equation}

Now the flow $g_t$ contracts vertical holonomy by $e^{-t}$ and expands horizontal holonomy by $e^t$.  By a slight abuse of notation we denote by $g_t(\gamma_n)$, the saddle connection $\gamma_n$ with respect to the Abelian differential $g_t(\omega_0)$.

These times are chosen so that for any $t\in (s_n,t_n)$ 
$$v(g_t(\gamma_n))=e^{-t}v(\gamma_n)\leq \frac{1}{\log n}$$ and 
$$\lambda(g_t(\gamma_n))=e^t\lambda(\gamma_n)\leq \frac{1}{\log n},$$
which implies that the length of $g_t(\gamma_n)$ is bounded by $\frac{2}{\log n}$ for any $t$ in this interval. 

If we can show that the intervals $(s_n,t_n)$ and $(s_{n+1},t_{n+1})$ overlap, then we would be done, since this would imply that at all times there is a short saddle connection. 
By (\ref{eq:saddle}), (\ref{eq:s}),(\ref{eq:t}) and  
the fact that $v(\gamma_n)\leq h_n(4)$ implies that  
\begin{equation}
\label{eq:time:interval}
t_n-s_n\geq 3\log n + O(\log\log n) 
\end{equation}

Now the form of the leading terms   of the matrix $C_{n+1}C_{n+2}$
implies that 
$$\frac{\lambda_n(3)}{\lambda_n(4)}\asymp O(n).$$
(Here, $A\asymp O(n)$ means $A=O(n)$ and $\frac{1}{A}=O(\frac{1}{n})$.)  
Using (\ref{eq:three:wins}) and the form of the leading terms 
of the matrix $C_{n+1}$ we have 
$$\frac{\lambda_n(3)}{\lambda_{n+1}(3)}= n^2+O(1/n)$$ so that 
\begin{equation}
\label{eq:successive}
\frac{\lambda_{n+1}(4)}{\lambda_n(4)}\asymp O(1/n^2).
\end{equation}

Now inequality~(\ref{eq:successive})
implies that $$t_{n+1}-t_n= 2\log n + O(1)$$ 
and together with (\ref{eq:time:interval}) implies that for large $n$, $$(s_n,t_n)\cap (s_{n+1},t_{n+1})\neq\emptyset.$$
This finishes the proof. 

\end{document}